\documentclass[12pt]{amsart}
\usepackage{amssymb,amscd}

\newtheorem{theorem}{Theorem}
\newtheorem{proposition}[theorem]{Proposition}

\newtheorem{lemma}[theorem]{Lemma}
\newtheorem{corollary}[theorem]{Corollary}

\theoremstyle{remark}
\newtheorem{remark}[theorem]{Remark}

\theoremstyle{definition}
\newtheorem{definition}[theorem]{Definition}

\newtheorem{convention}[theorem]{Convention}
\newtheorem{problem}[theorem]{Problem}

\numberwithin{equation}{section}
\numberwithin{theorem}{section}

\newcommand{\B}{\mathcal{B}}

\newcommand{\cone}{\widehat{\z(\M)}_+}

\newcommand{\h}{\mathfrak{H}}

\newcommand{\K}{\mathcal{K}}

\newcommand{\M}{\mathcal{M}}

\newcommand{\N}{\mathbb{N}}
\newcommand{\ozm}{\Omega(\z(\M))}

\newcommand{\p}{\mathcal{P}}
\newcommand{\ps}{\mbox{$(\p(\M)/\sim)$}}

\newcommand{\R}{\mathbb{R}}

\newcommand{\z}{\mathcal{Z}}

\begin{document}

\title[Dimension theory]{On the dimension theory of von Neumann algebras}
\author{David Sherman}
\address{Department of Mathematics\\ University of California\\ Santa Barbara, CA 93106}
\email{dsherman@math.ucsb.edu}

\begin{abstract}
In this paper we study three aspects of $\ps$, the set of Murray-von Neumann equivalence classes of projections in a von Neumann algebra $\M$.  First we determine the topological structure that $\ps$ inherits from the operator topologies on $\M$.  Then we show that there is a version of the center-valued trace which extends the dimension function, even when $\M$ is not $\sigma$-finite.  Finally we prove that $\ps$ is a complete lattice, a fact which has an interesting reformulation in terms of representations.
\end{abstract}

\maketitle

\section{Introduction} \label{S:intro}

Let $\M$ be a von Neumann algebra, $\p(\M)$ its projections, and $\sim$ the relation of Murray-von Neumann equivalence on $\p(\M)$.  The description of the quotient $\ps$ is known as the \textit{dimension theory} for $\M$.  In this paper we prove basic results about three aspects of dimension theory: topology, parameterization, and order.

The second section of the paper contains background which is relevant for all three topics.  Section 3 deals with topology; Sections 4 and 5 with parameterization; Sections 6 and 7 with order structure.  Except for one or two references, these three groupings are independent from each other.  In the remainder of this introduction we explain the problems which motivate our investigations.

\smallskip

\textsc{Topology.}  The first goal requires little explanation.
\begin{problem}
Study the topology that $\ps$ inherits from the strong (equivalently, the weak) topology on $\M$.
\end{problem}
\noindent Some of the results are used in the author's recent work on unitary orbits (\cite{S}).

\smallskip

\textsc{Parameterization.}  It is easy to check that $\ps$ also inherits a well-defined partial order from $\p(\M)$.  Classical work of Murray and von Neumann (\cite{MvN1}) and Dixmier (\cite{Di1949,Di1952}) shows that $\ps$ can be naturally parameterized by a subset of the extended positive cone of the center, at least when $\M$ is $\sigma$-finite.  This parameterization map, called a \textit{dimension function}, can be extended to all of $\M_+$, and the extension is called an \textit{extended center-valued trace}.  The existence of a dimension function on a non-$\sigma$-finite von Neumann algebra is also classical, though less-known.  It was originally studied in connection with spatial isomorphisms by Griffin (\cite{G1953,G1955}) and Pallu de la Barri\`{e}re (\cite{P}), and eventually given a representation-free foundation by Tomiyama (\cite{To}).

There is a noticeable gap between the last two objects.
\begin{problem}
Is there a version of the extended center-valued trace which extends the dimension function on a non-$\sigma$-finite von Neumann algebra?
\end{problem}
\noindent One might expect (and dread) technical constructions involving cardinals and limits.  We show how to avoid most of this by simply marrying Tomiyama's dimension function to the equivalence relation of Kadison and Pedersen (\cite{KP}).  In fact, the main point to settle does not involve cardinals.

\smallskip

\textsc{Order.}  The range of Tomiyama's map consists of certain cardinal-valued order-continuous functions on the spectrum of the center.  Tomi-yama assumed pointwise order and arithmetic on the range, then gave some examples to show that his map lacks basic continuity properties.  In fact the pointwise operations (on infinite sets of functions) do not behave well, and it seems to us that these are essentially the wrong operations to be considering.  Our viewpoint here is more algebraic.  This repairs certain degeneracies and allows us to resolve affirmatively the basic
\begin{problem} \label{P:complat}
Is $\ps$ always a complete lattice?
\end{problem}

We recall that a \textit{lattice} (resp. \textit{complete lattice}) is a partially-ordered set in which one may take meets and joins of finitely (resp. arbitrarily) many elements.  $\p(\M)$ is a complete lattice, but it does not induce lattice operations on $\ps$: for example, $[p] \wedge [q]$ is not well-defined as $[p \wedge q]$.  Nonetheless the comparison theorem for projections readily implies that $\ps$ is a lattice.  And in a finite von Neumann algebra, the dimension function identifies $\ps$ with a complete sublattice of $\z(\M)_1^+$.  Problem \ref{P:complat} asks about the existence of meets and joins of arbitrarily large sets of equivalence classes coming from arbitrarily large von Neumann algebras.  Its answer has a somewhat surprising reformulation in terms of representations.

\section{Background} \label{S:back}

Let $\M$ be a von Neumann algebra of arbitrary type and cardinality.  We write $\z(\M)$ for its center, and we occasionally symbolize the strong and weak topologies by $s$ and $w$.  The central support of an operator is $c(\cdot)$.

We use the standard terminology and results from \cite[Section V.1]{T} for projections, including $p^\perp$ for $(1-p)$.  Besides $p \sim q$, we write $p \preccurlyeq q$ for subequivalence, and $p \prec q$ for $p \preccurlyeq q$ but not $p \sim q$. 
Notice that for pairwise orthogonal sets $\{p_\alpha\}, \{q_\alpha\} \subset \p(\M)$,
\begin{equation} \label{E:addeq}
p_\alpha \sim q_\alpha, \: \forall \alpha \Rightarrow \left(\sum p_\alpha\right) \sim \left(\sum q_\alpha\right),
\end{equation}
\begin{equation} \label{E:addeq2}
p_\alpha \preccurlyeq q_\alpha, \: \forall \alpha \Rightarrow \left(\sum p_\alpha\right) \preccurlyeq \left(\sum q_\alpha\right).
\end{equation}

Among the many adjectives which may be applied to a single projection, we specify one which may cause confusion.  A nonzero projection $p$ is \textit{properly infinite} if $zp$ is infinite or zero for any central projection $z$.  (An alternative definition: $p$ is properly infinite if it can be decomposed into a countably infinite sum of projections, each of which is equivalent to $p$.)  Any adjective can be applied to an algebra when the adjective describes the identity projection of the algebra.

According to \eqref{E:addeq}, we can sum unambiguously any set in $\ps$ for which there are mutually orthogonal representatives, simply by taking the equivalence class of the sum of representatives.  This determines a partial order on $(\p(\M)/\sim)$: $[p] \leq [q]$ if there exists a projection $r$ with $[p] + [r] = [q]$.  One may also induce the same order directly, since the quotient operation respects the order in $\p(\M)$.  By this we mean
$$[p_1] \leq [p_2] \iff \exists q_1, q_2 \text{ with } q_1 \sim p_1, \: q_2 \sim p_2, \: q_1 \leq q_2.$$
So $[p_1] \leq [p_2]$ means nothing other than $p_1 \preccurlyeq p_2$.

Actually the comparison theorem for projections (\cite[Theorem V.1.8]{T}) implies that $(\p(\M)/\sim)$ is a lattice.  For $p,q \in \p(\M)$, let $z$ be a central projection with $zp \preccurlyeq zq$, $z^\perp p \succcurlyeq z^\perp q$.  Then 
\begin{equation} \label{E:lattice}
[p] \wedge [q] = [zp + z^\perp q], \qquad [p] \vee [q] = [z^\perp p + zq].
\end{equation}

Next we recall basic properties of the extended center-valued trace.  This material is due to Dixmier (\cite{Di1949,Di1952}), but for the reader's convenience (presumably), we give citations from Takesaki's book \cite{T}.

\begin{definition} (\cite[Definition V.2.33]{T}) Let $\M$ be an arbitrary von Neumann algebra, and let $\ozm$ be the spectrum of the abelian $C^*$-algebra $\z(\M)$.  By $\cone$ we mean the partially-ordered monoid of $[0,+\infty]$-valued continuous functions on $\ozm$.  $\z(\M)_+$ is contained in $\cone$ and acts on it by multiplication.

An \textbf{extended center-valued trace} on $\M$ is an additive map $T: \M_+ \to \cone$ which commutes with the action of $\z(\M)_+$ and satisfies $T(x^*x) = T(xx^*)$ for $x \in \M_+$.

$T$ is \textit{faithful} if $T(x^*x)=0 \Rightarrow x=0, \: \forall x \in \M_+$.  $T$ is \textit{normal} if
\begin{equation} \label{E:normal}
T(\sup x_\alpha) = \sup T(x_\alpha)
\end{equation}
for any bounded increasing net $\{x_\alpha\} \subset \M_+$.  $T$ is \textit{semifinite} if $\{x \in \M \mid T(x^*x) \in \z(\M)_+\}$ is $\sigma$-weakly dense in $\M$.
\end{definition}

Here we wish to draw attention to a point which will be amplified in Sections \ref{S:cont} and \ref{S:comp}.  What is the meaning of the expression $\sup T(x_\alpha)$ in \eqref{E:normal}?  The pointwise supremum of an increasing family of $[0,+\infty]$-valued continuous functions on $\ozm$ may not be continuous, and some kind of algebraic supremum is required instead.  Dixmier showed that such a supremum exists, using the fact that $\ozm$ is stonean (\cite{Di1951}).  He also mentions specifically that other methods, including a purely formal one, could reach the same goal (\cite[p.25]{Di1952}).  We suppose that our technique in Section \ref{S:comp} is similar to the formal approach that he had in mind. 

Semifinite von Neumann algebras - those with no summand of type III - are characterized by the existence of a faithful normal semifinite extended center-valued trace (\cite[Theorem V.2.34]{T}).  Such a map $T$ is unique up to multiplication by an element of $\cone$ which takes finite values on an open dense subset of $\ozm$, so all are equally useful in calculations.  A projection $p$ is finite if and only if $T(p)$ takes finite values on an open dense subset of $\ozm$ (\cite[Proposition V.2.35]{T}).  From all this $p \preccurlyeq q \Rightarrow T(p) \leq T(q)$, and the converse holds if $p$ is finite.

If $\M$ is finite, there is a \textit{unique} faithful extended center-valued trace $T$ with $T(1_\M) = 1_\M$ (\cite[Theorem V.2.6]{T}).  Such a map is automatically normal, and the linear extension which is defined on all of $\M$ is called simply a \textit{center-valued trace}.

\begin{convention} \label{C}
Whenever we talk of an  ``extended center-valued trace" $T$ on $\M_+$ in the sequel, it is assumed that
\begin{itemize}
\item $T$ is normal and faithful;
\item on the finite summand of $\M$, $T$ agrees with the center-valued trace;
\item on the semifinite summand of $\M$, $T$ is semifinite;
\item on the infinite type I summand of $\M$, $T$ maps an abelian projection to its central support.
\end{itemize}
Therefore $T(p) = (+\infty) c(p)$ for a projection supported on the type III summand.
\end{convention}

A word about operator topologies on $\M$: the strong, $\sigma$-strong, weak, and $\sigma$-weak topologies can all be defined spatially.  The $\sigma$-strong and $\sigma$-weak topologies are independent of the choice of (faithful normal) representation, and this is not true for the strong and the weak.  But on \textit{bounded} sets, we have the agreements strong=$\sigma$-strong and weak=$\sigma$-weak; we therefore permit ourselves the small linguistic abuse of referring to the strong (or weak) topology on a bounded subset of $\M$.

For $\M$ finite, the normality of the center-valued trace is equivalent to $\sigma$-weak-$\sigma$-weak continuity.  It will be more useful for us that this map is also $\sigma$-strong-$\sigma$-strong continuous, and therefore strong-strong continuous on bounded sets.  (See \cite[Theorem 13]{G1953}, \cite[I.4.Th\'{e}or\`{e}me 2 and p.250]{Di1969}, or \cite{R} in connection with this.  In fact the strong-strong or weak-weak continuity on all of $\M$ does depend on the representation (\cite[Theorem 8]{G1953}).)

\bigskip

Here are some examples of $\ps$.

\begin{enumerate}
\item When $\M$ is a type $\text{I}_n$ factor, $(\p(\M)/\sim)$ is isomorphic to the initial segment of cardinals $\leq n$, via the map that sends a projection to its rank.
\item When $\M$ is a type $\text{II}_1$ factor, $(\p(\M)/\sim) \simeq [0,1]$.
\item When $\M$ is a $\sigma$-finite type $\text{II}_\infty$ factor, $(\p(\M)/\sim) \simeq [0, +\infty]$.
\item When $\M$ is a $\sigma$-finite type $\text{III}$ factor, $(\p(\M)/\sim) \simeq \{0, +\infty\}$.
\end{enumerate}
The isomorphisms in (2) and (3) are implemented by a (bounded or unbounded) trace.  When $\M$ is a non-factor with separable predual, $(\p(\M)/\sim)$ is naturally viewed as a direct integral of the lattices above.  When $\M$ is finite, $(\p(\M)/\sim)$ is isomorphic to a sublattice of $\z(\M)_1^+$ via the center-valued trace (see \cite[Theorem 8.4.4]{KR2}).

Continuous (type II) and degenerate (type III) dimension theory were part of the original appeal for Murray and von Neumann: what happens at large cardinality?  Since $\ps$ is totally ordered if and only if $\M$ is a factor, this is the scenario closest to set theory.  Do type II and III factors contain ``quantum cardinal arithmetic" which diverges from the usual cardinal arithmetic of a type I factor?  

\smallskip

The questions above are answered neatly by the parameterization of $\ps$ as developed by Griffin (\cite{G1953,G1955}), Pallu de la Barri\`{e}re (\cite{P}), and especially as formulated by Tomiyama (\cite{To}).  The main point is a structure theorem allowing us to break a properly infinite von Neumann algebra into direct summands, each of which has a well-defined size.  This is in direct analogy to the structure theorem for type I von Neumann algebras, but we use $\sigma$-finiteness instead of abelianness as the ``unit of measurement".

\begin{definition} \label{D:homog} $\text{(\cite[Definition 1]{To})}$
Let $\kappa$ be a cardinal.  We say that a nonzero projection $p$ in a von Neumann algebra $\M$ is \textbf{$\kappa$-homogeneous} if $p$ is the sum of $\kappa$ mutually equivalent projections, each of which is the sum of centrally orthogonal $\sigma$-finite projections.  We also define 
$$\kappa_\M = \sup \{\kappa \mid \text{$\M$ contains a $\kappa$-homogeneous projection}\}.$$
\end{definition}

\begin{remark} The terminology here is conflicting.  We follow Tomiyama, but elsewhere ``$\kappa$-homogeneous projection" means a central projection which is the sum of $\kappa$ equivalent abelian projections (e.g. \cite[p.299]{T}).
\end{remark}

A projection can be $\kappa$-homogeneous for at most one $\kappa \geq \aleph_0$; also for $\kappa \geq \aleph_0$, two $\kappa$-homogeneous projections with identical central support are necessarily equivalent (\cite{G1955,To}).  $\kappa_\M$ is not larger than the dimension of a Hilbert space on which $\M$ is faithfully represented.

The fundamental result for us is a m\'{e}lange of two theorems of Griffin, one covering the semifinite case (slightly adapted to our setting, and also proved by Pallu de la Barri\`{e}re) and one covering the purely infinite.  It was rewritten in the non-spatial setting by Tomiyama.

\begin{theorem} \label{T:griffin} $($\cite[Theorem 3]{G1953}, \cite[Theorem 1]{G1955}, also \cite[I.5]{P} and \cite[Theorem 1]{To}$)$
Let $\M$ be a properly infinite von Neumann algebra.  Then uniquely
$$1_\M = \sum_{\aleph_0 \leq \kappa \leq \kappa_\M} z_\kappa,$$
where each $z_\kappa$ is either zero or a $\kappa$-homogeneous central projection.
\end{theorem}

Let $T$ be an extended center-valued trace on a von Neumann algebra $\M$ (following Convention \ref{C}).  Given any $p \in \p(\M)$, let $z^f$ be the largest central projection such that $z^f p$ is finite.  By applying Theorem \ref{T:griffin} to $(1-z^f) p\M p$, there are unique central projections $(z_\kappa)_{\aleph_0 \leq \kappa \leq \kappa_\M}$ such that $\sum z_\kappa p = (1-z^f)p$ and any nonzero $z_\kappa p $ is $\kappa$-homogeneous.  Make the formal assignment
\begin{equation} \label{E:gdf}
p = \left( z^f p + \sum_{\aleph_0 \leq \kappa \leq \kappa_\M} z_\kappa p \right) \mapsto \left( T(z^f p) + \sum_{\aleph_0 \leq \kappa \leq \kappa_\M} \kappa z_\kappa \right).
\end{equation}
From our earlier comments this assignment is a complete invariant for the equivalence class of $p$.

Under the isomorphism $\z(\M) \simeq C(\ozm)$, projections correspond to clopen subsets of $\ozm$, so elements on the right-hand side of \eqref{E:gdf} can be interpreted as partially-defined functions on $\ozm$.  The range is in $([0,+\infty) \cup \{\kappa \mid \aleph_0 \leq \kappa \leq \kappa_\M\})$, and the functions are (order) continuous on their domains, which are easily shown to be open and dense.  Tomiyama showed (\cite[Lemma 5]{To}) that such functions extend uniquely to continuous functions on all of $\ozm$.  

\begin{definition} \label{T:defgdf} (\cite{To})
The assignment described above, from $\p(\M)$ to the continuous $([0,+\infty) \cup \{\kappa \mid \aleph_0 \leq \kappa \leq \kappa_\M\})$-valued functions on $\ozm$, is a \textbf{(generalized) dimension function} of $\M$.
\end{definition}

\begin{theorem} \label{T:gdf} $($\cite{To}$)$
Let $D$ be a dimension function of $\M$.  Then $D$ is additive on pairs of orthogonal projections, provided that one incorporates the positive reals into cardinal arithmetic in the obvious way.  We have
$$p \preccurlyeq q \iff D(p) \leq D(q), \qquad \forall p,q \in \p(\M),$$
where we use the pointwise ordering of functions on the right-hand side.
\end{theorem}

It follows that $D$ factors as
$$\p(\M) \twoheadrightarrow \ps \overset{\sim}{\to} D(\p(\M)).$$
Here the second map is an embedding in a function space, preserving order, sums (when they exist), and the multiplicative $\p(\z(\M))$-action. 

\begin{corollary} \label{T:factor} ${}$
\begin{enumerate}
\item In a factor of type $\text{II}_\infty$, the totally ordered set $(\p(\M)/\sim)$ is isomorphic to
$$[0, +\infty) \cup \{\kappa \mid \aleph_0 \leq \kappa \leq \kappa_\M\}.$$
\item In a factor of type III, the totally ordered set $(\p(\M)/\sim)$ is isomorphic to
$$\{0\} \cup \{\kappa \mid \aleph_0 \leq \kappa \leq \kappa_\M\}.$$
\end{enumerate}
\end{corollary}

So any interest in ``quantum cardinal arithmetic" wanes here: infinite quantum cardinals are (isomorphically) just cardinals.  For the reader interested in axiomatic treatments of $\ps$ and more general algebraic structures obtained as quotients of lattices, see \cite{L,M,F}.

\section{The topology of $(\p(\M)/\sim)$} \label{S:top}

If we want $\ps$ to inherit a topology from $\p(\M)$, there really are not so many interesting choices.  The quotient of the norm topology is the discrete topology, since $\|p - q\| < 1$ implies that $p$ and $q$ are unitarily equivalent (\cite[5.2.6-10]{W-O}).  And all of the ``operator" topologies (notably, the strong and the weak) are equivalent when restricted to $\p(\M)$ (\cite[Ex. 5.7.4]{KR1}).  We point out, however, that $(\p(\M), \text{strong})$ is complete, while $(\p(\M), \text{weak})$ may not be; completeness is not a topological property.

We will denote the resulting quotient strong/weak operator topology on $(\p(\M)/\sim)$ by ``$QOT$".  In the rest of this section, all closures and convergences in $\ps$ are to be understood in this topology.

We need a few lemmas.

\begin{lemma} \label{T:liminf}
Let $\{x_\alpha\}$ be a net in a semifinite von Neumann algebra $\M$ equipped with an extended center-valued trace $T$.  If $x_\alpha^* x_\alpha = y_1$ is fixed, while $x_\alpha x_\alpha^* \overset{w}{\to} y_2$, then $T(y_1) \geq T(y_2)$ in $\widehat{\z(\M)}_+$.
\end{lemma}

\begin{proof}
Fix any $\varphi \in \z(\M)_*^+$.  Then $\varphi \circ T$ is a semifinite normal weight, so weakly lower-semicontinuous (\cite{H}).  We have
\begin{align*}
\varphi \circ T(y_2) &= \varphi \circ T(w-\lim x_\alpha x_\alpha^*) \leq \liminf \varphi \circ T(x_\alpha x_\alpha^*) \\ &= \liminf \varphi \circ T(x_\alpha^* x_\alpha) = \varphi \circ T(y_1).
\end{align*}
Since $\varphi$ is arbitrary, the conclusion follows.
\end{proof}

\begin{lemma} \label{T:proj}
Let $p,q,r \in \p(\M)$, with $\M$ and $p$ properly infinite.
\begin{enumerate}
\item If $p \sim q_j$ for a countable set $\{q_j\}$, then $p \sim \vee q_j$.
\item If $zq \prec zr$ for all nonzero central projections $z$, then $q^\perp \sim 1_\M$.
\end{enumerate}
\end{lemma}

\begin{proof} $\quad$

(1) It is clear that $p \preccurlyeq \vee q_j$.  Write $p = \sum p_j$, where each $p_j \sim p$.  Let $v_j$ be a partial isometry between $p_j$ and $q_j$.  The operator $\sum v_j/2^j$ has right support $\vee q_j$ and left support $\leq p$, so also $p \succcurlyeq \vee q_j$. 

(2) We compare $q$ and $q^\perp$.  If there were a nonzero central projection $z$ with $zq \succcurlyeq zq^\perp$, then $zq$ would be properly infinite (else a nonzero central projection would be the sum of two finite projections).  Write $zq = zq_1 + zq_2$, where $zq_1 \sim zq_2 \sim zq$.  By \eqref{E:addeq2},
$$zq \preccurlyeq z = (zq + zq^\perp) \preccurlyeq (zq_1 + zq_2) = zq,$$
so that $zq \sim z$.  Now $zr \succ zq \sim z$, which is impossible.

Thus $q \preccurlyeq q^\perp$.  By the same argument, $q^\perp$ is properly infinite and $q^\perp \sim 1_\M$.
\end{proof}

\begin{theorem} \label{T:sinfin}
If $\M$ is a finite von Neumann algebra, the center-valued trace induces a homeomorphism from $(\ps, QOT)$ to a subspace of $(\z(\M)_1^+, \text{\textnormal{strong}})$.  Consequently
$$\overline{\{[p]\}} = \{[p]\}, \qquad p \in \p(\M).$$
\end{theorem}

\begin{proof} Let $T$ be the center-valued trace.  If $[p_\alpha] \to [p]$, then there exist $q_\alpha \sim p_\alpha$ with $q_\alpha \overset{s}{\to} p$.  By the strong-strong continuity of $T$ on bounded sets, we have $T(p_\alpha) = T(q_\alpha) \overset{s}{\to} T(p)$.

On the other hand, suppose $p_\alpha, p$ are projections such that $T(p_\alpha) \overset{s}{\to} T(p)$.  Let $q_\alpha \leq p$ be projections with $T(q_\alpha) = T(p_\alpha) \wedge T(p)$, where the meet is taken in $\z(\M)_1^+$.  Let $r_\alpha \perp q_\alpha$ be projections with $T(r_\alpha) = T(p_\alpha - p) \vee 0$.  It follows that $r_\alpha$ is centrally orthogonal to $(p-q_\alpha)$, and by comparing center-valued traces $(q_\alpha + r_\alpha) \sim p_\alpha$.

When $\M$ is $\sigma$-finite, the strong topology on bounded sets is generated by the norm $x \mapsto \tau(x^*x)^{1/2}$, for $\tau$ any faithful tracial state (\cite[Proposition III.V.3]{T}).  A general finite algebra is a direct sum of $\sigma$-finite ones (\cite[Corollary V.2.9]{T}), so it suffices to show convergence for the seminorms coming from a family of traces, each of which is faithful on a $\sigma$-finite summand.

We now take such a trace $\tau$ and compute
\begin{align*}
\tau([(q_\alpha + r_\alpha) - p]^2) &= \tau([r_\alpha - (p - q_\alpha)]^2) = \tau (r_\alpha + (p-q_\alpha)) \\ &= \tau (T(r_\alpha) + T(p-q_\alpha)) = \tau(|T(r_\alpha) - T(p-q_\alpha)|) \\ &= \tau(|T((q_\alpha + r_\alpha) - p)|) = \tau(|T(p_\alpha - p)|) \\ &\leq \tau(|T(p_\alpha - p)|^2)^{1/2} \to 0. \qedhere
\end{align*}
\end{proof}

Regarding Theorem \ref{T:sinfin}, we remind the reader that typically we do \textit{not} have an equivalence between the strong and weak topologies on $\z(\M)_1^+$.

\begin{theorem} \label{T:sininf}
Let $p$ be a projection in a properly infinite von Neumann algebra $\M$.  If $p$ is finite, 
\begin{equation} \label{E:finclos}
\overline{\{[p]\}} = \{[q] \mid [q] \leq [p]\}.
\end{equation}
If $p$ is properly infinite and $c(p)=1_\M$,
\begin{equation} \label{E:infclos}
\overline{\{[p]\}} = (\p(\M)/\sim).
\end{equation}
Equations \eqref{E:finclos} and \eqref{E:infclos} may be synthesized into
\begin{equation} \label{E:seg}
\overline{\{[p]\}} = \{[q] \mid T(q) \leq T(p)\}, \qquad \forall p \in \p(\M),
\end{equation}
for any extended center-valued trace $T$.
\end{theorem}

\begin{proof} First consider a finite projection $p$.  We may assume that $c(p) = 1$ and so $\M$ is semifinite; let $T$ be an extended center-valued trace.  If $p_\alpha \sim p$ and $p_\alpha \overset{w}{\to} q$, then by Lemma \ref{T:liminf}, $T(q) \leq T(p)$.  We have assumed $p$ finite, so $q$ is as well and $p \succcurlyeq q$.  For the other containment, choose any $q$ with $[q] \leq [p]$.  Write $p = p_0 + p_1$, with $p_0 \sim q$.  Since $q$ is finite, $q^\perp \sim 1$ is properly infinite, and we may write $q^\perp = \sum_{k=1}^\infty q_k$, with $q_k \sim q^\perp \sim 1$.  Let $p_1 \sim r_k \leq q_k$.  Then $p = (p_0 + p_1) \sim (q + r_k) \overset{s}{\to} q$.  This proves \eqref{E:finclos}.

Now consider arbitrary $q$ and properly infinite $p$ with $c(p)=1$.  Find the largest central projection $z$ with $zp \preccurlyeq zq$.  Consider the nonempty net $\{zp_\alpha \mid zp \sim zp_\alpha \leq zq\}$, with order inherited from $\p(\M)$.  It is upward directed by Lemma \ref{T:proj}(1), applied to two projections.  Its supremum is $zq$.

By Lemma \ref{T:proj}(2) $z^\perp q^\perp \sim z^\perp$, which is properly infinite, so we may write $z^\perp q^\perp$ as the countable sum $\sum q_j$, with each $q_j \sim z^\perp q^\perp \sim z^\perp$.  Write $z^\perp p = z^\perp p_0 + z^\perp p_1$, where $z^\perp p_0 \sim z^\perp q$.  Also let $z^\perp p_1 \sim r_j \leq q_j$.  Then $z^\perp p = (z^\perp p_0 + z^\perp p_1) \sim (z^\perp q + r_k) \overset{w}{\to} z^\perp q$.
  
Combining the results for $zp$ and $z^\perp p$ and considering the product net, we see that $q$ is a strong limit of projections equivalent to $p$.  This proves \eqref{E:infclos}.

Equation \eqref{E:seg} follows from \eqref{E:finclos} and \eqref{E:infclos} by breaking off the largest central summand where $p$ is properly infinite with full central support.
\end{proof}

\begin{corollary} \label{T:Eclos}
Let $\M$ be a factor and $E \subseteq \ps$.  We consider an extended center-valued trace $T$ on $\M$ to be a $[0,+\infty]$-valued function.

If $\M$ is finite,
$$\overline{E} = \{[q] \mid T(q) \in \overline{\{T(p) \mid [p] \in E\}}\}.$$

If $\M$ is properly infinite,
$$\overline{E} = \{[q] \mid T(q) \leq \sup_{[p] \in E} T(p)\}.$$
\end{corollary}

Corollary \ref{T:Eclos} follows readily from the preceding arguments, and its easy proof is left to the interested reader.

\begin{corollary}
$QOT$ is a $T_1$ topology exactly when $\M$ is finite.
\end{corollary}

\begin{proof}
This is a direct consequence of Theorems \ref{T:sinfin} and \ref{T:sininf}. 

A topology is $T_1$ if for any two distinct points $x,y$, there is a closed set which contains $x$ and not $y$.  If $\M$ is not finite, let $x$ be the equivalence class of a properly infinite projection, and let $y$ be $[0]$.  Since $y$ belongs to the closure of $x$, no such separating closed set exists.  (In general, a topology is $T_1$ iff singletons are closed.)
\end{proof}

It turns out to be more useful for our applications elsewhere (\cite{S}) to know when $QOT$ is $T_0$.  A topology is $T_0$ if for any two distinct points, there exists a closed set which contains exactly one of them.

\begin{proposition} \label{T:T0}
For a von Neumann algebra $\M$, the following conditions are equivalent.
\begin{enumerate}
\item $QOT$ is a $T_0$ topology on $(\p(\M)/\sim)$.
\item For any $p,q \in \p(\M)$, $[p] \in \overline{\{[q]\}} \Rightarrow p \preccurlyeq q$.
\item $\kappa_\M \leq \aleph_0$.
\item $\M$ is a (possibly uncountable) direct sum of $\sigma$-finite von Neumann algebras.
\item $\M$ does not contain $\B(\h_1)$, where $\h_1$ is a Hilbert space of dimension $\aleph_1$.
\end{enumerate}
\end{proposition}

\begin{proof} The equivalence of conditions (3)-(5) follows from the definitions and Theorem \ref{T:griffin}.  We therefore focus on the equivalence of (1)-(3).

(1) $\to$ (3): If (3) fails, let $q$ be an $\aleph_1$-homogeneous projection, and let $p$ be an $\aleph_0$-homogeneous projection with $c(p)=c(q)$.  Then $[p] \in \overline{\{[q]\}}$ and $[q] \in \overline{\{[p]\}}$, but $[p] \neq [q]$.  Clearly there is no closed set separating the two.

(3) $\to$ (2): When $\kappa_\M \leq \aleph_0$, $T |_{\p(\M)}$ can be identified with $D$.  By Theorems \ref{T:sinfin} and \ref{T:sininf} we have
$$[p] \in \overline{\{[q]\}} \Rightarrow T(p) \leq T(q) \Rightarrow D(p) \leq D(q) \Rightarrow p \preccurlyeq q.$$

(2) $\to$ (1): Suppose (2) holds.  Given $[p],[q] \in \ps$, they can be separated by a closed set if $[p] \notin \overline{\{[q]\}}$ or $[q] \notin \overline{\{[p]\}}$.  If neither of these is true, then
$$[p] \in \overline{\{[q]\}}, \: [q] \in \overline{\{[p]\}} \: \Rightarrow p \preccurlyeq q, \: q \preccurlyeq p \: \Rightarrow [p]=[q]. \qedhere$$
\end{proof}

\section{From dimension function to trace in full generality} \label{S:trace}

Let $T$ be an extended center-valued trace on a von Neumann algebra $\M$, with $D$ the induced dimension function.  We will create a map which extends $D$ to the entire positive cone and so is a trace which distinguishes among infinite cardinalities.  (In case $\kappa_\M \leq \aleph_0$, this process simply recovers $T$.)  The main tool is

\begin{definition} (\cite{KP})
For two elements $h,k \in \M_+$, we write $h \approx k$ if and only if there exists a family $\{x_\alpha\} \subset \M$ such that $h = \sum x_\alpha^* x_\alpha$ and $k = \sum x_\alpha x_\alpha^*$.

We write $h \lessapprox k$ to mean that there exists $k' \leq k$ with $h \approx k'$.

For $h \in \M_+$, we say that $h$ is \textit{finite} if $h \approx k \leq h \Rightarrow k=h$. 
\end{definition}

The following facts are shown in \cite{KP}.

\begin{itemize}
\item The relation $\approx$ is an equivalence relation.  It is homogeneous ($h \approx k \Rightarrow \lambda h \approx \lambda k, \: \lambda \in \R_+$) and completely additive in the sense that
$$h_\alpha \approx k_\alpha, \: \forall \alpha \quad \Rightarrow \quad \sum h_\alpha \approx \sum k_\alpha$$
(when the two sums exist in $\M$).

\item The relation $\lessapprox$ gives a partial order on equivalence classes.  In particular,
\begin{equation} \label{E:pord}
h \lessapprox k, \: k \lessapprox h \: \Rightarrow \: h \approx k, \qquad h,k \in \M_+.
\end{equation}

\item For projections, $p \approx q \iff p \sim q$.

\item For $h,k \in \M_+$, $h \lessapprox k \Rightarrow T(h) \leq T(k)$, and the converse holds if $h$ is finite.

\end{itemize}

We will also say that nonzero $h \in \M_+$ is \textit{properly infinite} if  $zh$ is finite and nonzero for no central projection $z$.  For projections, the usage here of ``finite" and ``properly infinite" coincides with the usual meaning; in fact proper infiniteness of (nonzero) $h$ in either case is characterized by $T(h)$ being $\{0,+\infty\}$-valued.

\begin{lemma} \label{T:mult} ${}$
\begin{enumerate}
\item Let $\lambda \in ((0,1) \cup (1,\infty))$, and let $p$ be a projection.  Then
$$p \text{ is properly infinite } \iff p \approx \lambda p.$$
\item Let $h,k \in \M_+$ have equal central support, with $k$ properly infinite and $h$ a countable sum of finite elements.  Then $h \lessapprox k$.
\item Let $h,k \in \M_+$ be properly infinite with equal central support, and suppose that each is a countably infinite sum of finite elements.  Then $h \approx k$.
\end{enumerate}
\end{lemma}

\begin{proof} (1) If $p \approx \lambda p$, then $T(p)$ must be $\{0,+\infty\}$-valued.  For the opposite implication, we first check rational multiples.  Let $m,n \in \N$.  By proper infiniteness, we may write
$$p = \sum_{i=1}^m p_i = \sum_{j=1}^n p'_n, \qquad p_i \sim p \sim p'_j, \: \forall i,j.$$
Then
\begin{align*}
p &= \sum_{i=1}^m p_i \approx \sum_{i=1}^m p = mp = \left( \frac{m}{n} \right) np =  \left( \frac{m}{n} \right) \left(\sum_{j=1}^n p \right) \\ &\approx \left( \frac{m}{n} \right) \left( \sum_{j=1}^n p'_n \right) = \left( \frac{m}{n} \right) p.
\end{align*}
Find two positive rationals $\lambda_1, \lambda_2$ with $\lambda_1 \leq \lambda \leq \lambda_2$:
$$p \approx \lambda_1 p \leq \lambda p \leq \lambda_2 p \approx p \; \Rightarrow \; p \approx \lambda p,$$
using \eqref{E:pord}.

(2) Write $h = \sum_{j=1}^\infty h_j$, where each $h_j$ is finite.  Since $T(h_1) \leq T(k)$, there is an operator $k_1$ with $h_1 \approx k_1 \leq k$.  We continue in this way: since $T(h_n) \leq T(k - \sum_{j=1}^{n-1} k_j)$, find $k_n$ with $h_n \approx k_n \leq (h - \sum_{j=1}^{n-1} k_j)$.

Now each $(\sum_{j=1}^n h_j) \approx (\sum_{j=1}^n k_j)$, and these terms are finite and increasing to $h$ and some $k'$, respectively.  It follows from \cite[Lemma 3.3]{KP} that $h \approx k' \leq k$.

(3) Both $h \lessapprox k$ and $h \gtrapprox k$ follow from the previous part; apply \eqref{E:pord}.
\end{proof}

\begin{proposition} \label{T:appproj}
Let $h \in \M_+$ be properly infinite.  Then there exists $p \in \p(\M)$ such that $h \approx p$.
\end{proposition}

\begin{proof} It does no harm to assume that $h$ has full central support, and therefore $\M$ is properly infinite.  Write the identity as $1_\M = \sum_{n=-\infty}^\infty p_n, \: 1_\M \sim p_n$, and let $r_0 \leq p_0$ be an $\aleph_0$-homogeneous projection with full central support.

Now make the decomposition
$$h = \sum_{n=1}^\infty (2^{-n} \|h\|) q_n,$$
where $q_n$ is the spectral projection for $h$ corresponding to
$$\bigcup_{j=1}^{2^{n-1}} \left( (2j-1)2^{-n}\|h\|, (2j) 2^{-n} \|h\| \right].$$

For each $n \geq 1$, let $z^f_n$ be the largest central projection such that $z^f_n q_n$ is finite.  Using Lemma \ref{T:mult}(1) and then conjugating by a partial isometry from $(1 - z^f_n)$ to $(1-z^f_n)p_n$, find a projection $r_n$ with
$$(1-z_n^f)(2^{-n} \|h\|) q_n \approx (1-z^f_n) q_n \sim r_n \leq p_n.$$
Conjugating by a partial isometry from $z^f_n$ to $z^f_n p_{-n}$, let $r_{-n}$ be any operator (necessarily finite, but not necessarily a projection) with
$$z^f_n (2^{-n} \|h\|) q_n \approx r_{-n} \in p_{-n} \M p_{-n}.$$
By construction we have $h \approx \sum_{n=1}^\infty (r_n + r_{-n})$.

Set $z_0 = \wedge z^f_n$.  We will complete the proof by showing that $z_0 h$ and $z_0^\perp h$ are both (Kadison-Pedersen) equivalent to projections.

First,
$$z_0 h \approx z_0 \left( \sum_{n=1}^\infty (r_n + r_{-n}) \right) = z_0 \left( \sum_{n=1}^\infty r_{-n} \right).$$
The left-hand side has central support $z_0$, and is either zero or properly infinite because $h$ is properly infinite.  The right-hand side is a countable sum of finite elements.  By Lemma \ref{T:mult}(3),
$$z_0 h \approx z_0 r_0.$$

Second,
$$z_0^\perp \left( \sum_{n=1}^\infty r_n \right) \sim z_0^\perp \left( r_0 + \sum_{n=1}^\infty r_n \right),$$
since the central supports are equal and the left-hand side is a properly infinite projection.  (For example, this follows by evaluating the dimension function on both sides and noting that adding $\aleph_0$ does not change an infinite cardinal.)  On the other hand, Lemma \ref{T:mult}(2) implies
$$z_0^\perp \left( \sum_{n=1}^\infty r_{-n} \right) \lessapprox z_0^\perp r_0.$$
We put these together:
\begin{align*}
z_0^\perp h &\approx z_0^\perp \left( \sum_{n=1}^\infty (r_n + r_{-n}) \right) \lessapprox z_0^\perp \left( r_0 + \sum_{n=1}^\infty r_n \right) \\ &\sim z_0^\perp \left(\sum_{n=1}^\infty r_n \right) \lessapprox z_0^\perp \left( \sum_{n=1}^\infty (r_n + r_{-n}) \right) \approx z_0^\perp h.
\end{align*}
Then all terms above are (Kadison-Pedersen) equivalent, and the middle two are projections.
\end{proof}

\begin{corollary} \label{T:absorb}
Under the same hypotheses as in Lemma \ref{T:mult}(2), $k \approx \lambda k$ for any $\lambda \in (0, \infty)$, and $(h + k) \approx k$.
\end{corollary}

\begin{proof}
By Proposition \ref{T:appproj} and Lemma \ref{T:mult}(1), there is a properly infinite projection $p$ with $k \approx p \approx \lambda p \approx \lambda k$.  By Lemma \ref{T:mult}(2),
$$(h + k) \lessapprox 2k \approx k \lessapprox (h+k) \: \Rightarrow \: (h+k) \approx k.$$
\end{proof}

We are now ready to define our map.

\begin{definition}
With $T$ (and $D$) given, we construct a \textbf{fully extended center-valued trace} $\widehat{T}$ on $\M$ as follows.

For any $h \in \M_+$, let $z^f$ be the largest central projection so that $z^f h$ is finite.  Let $p$ be a projection with $p \approx (1 -z^f) h$.  Such a $p$ exists by Proposition \ref{T:appproj}, and all choices belong to the same Murray-von Neumann equivalence class.

We define
\begin{equation}
\widehat{T}(h) = T(z^f h) + D((1- z^f) p),
\end{equation}
which we view as a continuous $([0,+\infty) \cup \{\kappa \mid \aleph_0 \leq \kappa \leq \kappa_\M\})$-valued function on $\ozm$.
\end{definition}

\begin{theorem} \label{T:Tfull}
The map $\widehat{T}$ extends $D$, is additive, commutes with the multiplicative action of $\z(\M)_+$, and satisfies
\begin{equation} \label{E:pres}
h \lessapprox k \iff \widehat{T}(h) \leq \widehat{T}(k), \qquad h,k \in \M_+.
\end{equation}
(We are allowing cardinal arithmetic to incorporate the positive reals in the obvious way.)
\end{theorem}

\begin{proof}
Clearly $\widehat{T}$ extends $D$.  By the properties of $D$ and $T$ we have $h \approx k \iff \widehat{T}(h) = \widehat{T}(k)$.

In saying that $\widehat{T}$ is additive, we mean that
\begin{equation} \label{E:add}
\widehat{T}(h + k) = \widehat{T}(h) + \widehat{T}(k), \qquad h,k \in \M_+.
\end{equation}
For $h,k$ finite, \eqref{E:add} follows from additivity of $T$.  For $h,k$ properly infinite, the projection representing $h + k$ may be constructed as the sum of orthogonal representing projections for $h$ and $k$; \eqref{E:add} then follows from the additivity of $D$.  Finally, let $h$ and $k$ have the same central support, with $h$ finite and $k$ properly infinite.  In this case $\widehat{T}(h)$ is bounded above by $\aleph_0$, while $\widehat{T}(k) \geq \aleph_0$ where it is nonzero.  So $\widehat{T}(h)  + \widehat{T}(k) = \widehat{T}(k)$.  Since $(h + k) \approx k$ by Corollary \ref{T:absorb}, $\widehat{T}(h + k) = \widehat{T}(k)$ as well.

In saying that $\widehat{T}$ commutes with the action of $\z(\M)_+$, we mean
\begin{equation} \label{E:modmap}
y \widehat{T}(h) = \widehat{T}(yh), \qquad y \in \z(\M)_+, \: h \in \M_+.
\end{equation}
Clearly \eqref{E:modmap} holds for finite elements, since the analogous formula is true for $T$.  It therefore suffices to prove \eqref{E:modmap} under the assumption that $h$ and $y$ have full central support, with $h$ properly infinite.  In this case $y \widehat{T}(h) = \widehat{T}(h)$, so we are left to show that $yh \approx h$.  If $y \geq \lambda c(y)$ for some $\lambda > 0$, then by Corollary \ref{T:absorb}
$$h \approx \lambda h \leq yh \leq \|y\| h \approx h \: \Rightarrow \: h \approx yh.$$
The general conclusion follows by writing $y$ as a central sum of operators which are invertible on their supports.

As for \eqref{E:pres}, the forward implication is a consequence of additivity.  For the reverse implication, we look at central summands: where $h$ is finite, this is a property of $T$; where $h$ and $k$ are both infinite, this is a property of $D$.
\end{proof}

From Theorem \ref{T:Tfull}, we see that $\widehat{T}$ factors as
$$\M_+ \twoheadrightarrow (\M_+/\approx) \overset{\sim}{\to} \widehat{T}(\M_+).$$
Here the second map is an embedding in a function space, preserving order, sums, and the multiplicative $\z(\M)_+$-action.

More generally, we may say that an arbitrary completely additive map on $\M_+$ which respects the $\R_+$-action is \textit{tracial} if and only if it factors through the quotient $\M_+ \twoheadrightarrow (\M_+/\approx)$.  Numerical (completely additive) traces result when the range is $[0,+\infty]$; they are ``one-dimensional representations" of $(\M_+/\approx)$.

\begin{remark} \label{R:trace} Kadison and Pedersen observed that all extended center-valued traces on semifinite algebras can be generated in the following manner (\cite[Theorem 3.8]{KP}).  Fix a finite projection $p$ with full central support, and assume that $p$ is the identity on the finite summand and abelian on the infinite type I summand (to match Convention \ref{C}).  Then for finite $h \in \M_+$, $T(h)$ is the unique element of the extended center with $h \approx T(h)p$.  Already this requires a small extension of $\approx$ to unbounded sums.

With a further extension involving cardinals, $\widehat{T}$ can also be defined in this way.  For general $\M$, let $p$ be the identity on the finite summand, abelian on the infinite type I summand, finite on the type II summand, and $\aleph_0$-homogeneous on the type III summand; of course $p$ should have full central support.  For $h \in \M_+$, one can define $\widehat{T}(h)$ as the unique formal sum (as in \eqref{E:gdf}) such that $h \approx \widehat{T}(h)p$ and $\widehat{T}(h)$ takes no finite nonzero values on the type III summand.  Probably this is more interesting to mention than to carry out, so we omit the details.
\end{remark}

\section{Continuity} \label{S:cont}

In the remaininder of the paper we assume that $T$, $D$, and $\widehat{T}$ are given on $\M$.

The order-preserving embeddings of $\ps$ and $(\M_+/\approx)$ in a function space (albeit cardinal-valued) make pointwise operations available.  From Theorems \ref{T:gdf} and \ref{T:Tfull} we know that for finite sets, addition in the quotient structures agrees with addition of functions.  One may likewise add up infinite sets of functions, but there is no guarantee that the sum will be continuous.  Tomiyama gave an example (\cite[Example 2]{To}) to show that for a pairwise orthogonal set $\{p_\alpha\}$, one cannot expect an identity between $\sum D(p_\alpha)$ and $D(\sum p_\alpha)$, so that $D$ is not completely additive.

This is really an artifact of the function representation.  $\ps$ carries a natural (partially-defined) sum operation, given by
$$\sum [p_\alpha] \triangleq \left[ \sum q_\alpha \right]$$
whenever there exists a set of pairwise orthogonal projections $\{q_\alpha\}$ with $q_\alpha \sim p_\alpha$.  A similar definition is possible for sums in $\widehat{T}(\M_+)$, where we simply require that the representatives sum to an element of $\M_+$.  Note that there is no ambiguity in these definitions, by \eqref{E:addeq} and the definition of $\approx$, and as an immediate consequence, the maps $\p(\M) \twoheadrightarrow \ps$ and $\M_+ \twoheadrightarrow (\M_+/\approx)$ are completely additive.  It is of course possible to transport these sum operations to $D(\p(\M))$ and $\widehat{T}(\M_+)$.   

\smallskip

Pointwise lattice operations on pairs in $D(\p(\M))$ match \eqref{E:lattice} and so agree with the operations in $\ps$, but meets and joins of infinite sets of continuous functions need not be continuous.  For bounded real-valued functions on a stonean space, a regularization corrects this problem (\cite[Section III.1]{T}), but the situation for cardinal-valued functions is less clear.

Normality for $D$ and $\widehat{T}$ means an appropriate analogue of \eqref{E:normal}.  So how do we interpret an expression like $\sup D(p_\alpha)$, where $\{p_\alpha\}$ is an increasing net in $\p(\M)$?  As we just mentioned, the pointwise supremum need not lie in $D(\p(\M))$.  In the next section we show that the supremum always does make sense in $\ps$, but unfortunately normality is to much to ask.  Tomiyama gave a simple example (\cite[Example 1]{To}) of an uncountable increasing family of projections $\{p_\alpha\}$ for which the pointwise supremum of $D(p_\alpha)$ does lie in $D(\p(\M))$, and yet $\sup D(p_\alpha) \neq D(\sup p_\alpha)$.

This represents a phenomenon which really occurs in the quotient map $\p(\M) \twoheadrightarrow \ps$, as we already saw in the proof of \eqref{E:infclos}.  For $p$ a properly infinite projection, the elements of $[p]$, under the operator ordering, form an increasing net which converges strongly to $c(p)$.  One obtains a counterexample to normality whenever $c(p) \notin [p]$, and such counterexamples exist when $\kappa_\M > \aleph_0$.  On the other hand, if $\kappa_\M \leq \aleph_0$, the quotient maps are given by the extended center-valued trace, which we know to be normal.  We conclude

\begin{proposition} \label{T:normal}
Another equivalent condition in Proposition \ref{T:T0} is
\begin{enumerate}
\item[(6)] The quotient maps $\p(\M) \twoheadrightarrow \ps$ and $\M_+ \twoheadrightarrow (\M_+/\approx~)$ are normal.
\end{enumerate}
\end{proposition}

In contrast, a pointwise criterion for normality of $D$ and $\widehat{T}$ holds if and only if $\kappa_\M \leq \aleph_0$ and the center of $\M$ is finite-dimensional.  We do not bother to prove this explicitly, but we mention an example.  Let $\M = \ell^\infty$, and take $p_n$ to be the sum of the first $n$ elements of the standard basis.  Since $\sup D(p_n)$ does not agree with $D(\sup p_n)$ at any point of $(\beta \N \setminus \N) \subset \beta \N \simeq \ozm$, pointwise normality fails.  And here $D$ is the identity!

\smallskip

Our conclusion from all this is that the pointwise lattice and addition operations on functions in the range of $D$ and $\widehat{T}$ should be shelved in favor of the induced quotient structures on $\ps$ and $(\M_+/\approx)$.  With this interpretation the assertion ``$D$ and $\widehat{T}$ are normal" is also equivalent to the conditions in Proposition \ref{T:T0}.

\section{$\ps$ is a complete lattice} \label{S:comp}
Having just been warned about the degeneracies of the pointwise ordering, we omit the last step of Tomiyama's construction for $D$ and stick with a more algebraic language.  We follow the right-hand side of \eqref{E:gdf}, further dividing $T(z^f p)$ into the pieces where it lies between consecutive finite cardinals.  This allows us to write the typical element of $\widehat{T}(\M_+)$ as
\begin{equation} \label{E:simpler}
\sum_{\kappa \leq \kappa_\M} g_\kappa z_\kappa.
\end{equation}
The meaning of this expression is as follows.  If $\kappa$ is an infinite cardinal, then $g_\kappa = \kappa$.  If $\kappa$ is a nonnegative integer, $g_\kappa$ is an element of $\z(\M)_+$ satisfying $(\kappa-1)z_\kappa \leq g_\kappa \leq \kappa z_\kappa$ and $c(g_\kappa - (\kappa-1)z_\kappa) = z_\kappa$.  The central projections $z_\kappa$ sum to 1, and the decomposition is unique.

The partial order, pairwise sum operation, and pairwise lattice operations are all easily implemented for expressions of the form \eqref{E:simpler}.  Conversely, such an expression belongs to $\widehat{T}(\M_+)$ if it is $\{0, +\infty\}$-valued on the type III summand and less than some finite multiple of $\widehat{T}(1_\M)$.  To belong to $D(\p(\M)) \subseteq \widehat{T}(\M_+)$, an expression must be $\leq D(1_\M)$ and appropriately valued on both the type I and type III summands.

\begin{theorem} \label{T:complat}
$(\M_+/\approx)$ and $\ps$ are complete lattices.
\end{theorem}

\begin{proof}
We show how to perform lattice operations on formal sums of the form \eqref{E:simpler}.  Our constructions will preserve all conditions mentioned in the paragraph before Theorem \ref{T:complat}, so they are well-defined operations in $(\M_+/\approx)$ and $\ps$.

Let us find the supremum of an arbitrary set $\{f^\alpha\}$, where 
$$f^\alpha = \sum g_\kappa^\alpha z_\kappa^\alpha.$$

For each cardinal $\kappa \leq \kappa_\M$, set
$$y_{\leq \kappa} = \bigwedge_\alpha \left(\sum_{\lambda \leq \kappa} z^\alpha_\lambda \right);$$
$y_{\leq \kappa}$ is ``where all $f^\alpha$ are $\leq \kappa$".  Note that $y_{\leq \kappa}$ is increasing in $\kappa$ and $y_{\leq \kappa_\M} = 1$.  Next define, for each cardinal $\kappa \leq \kappa_\M$,
$$z_\kappa = y_{\leq \kappa} - \bigvee_{\lambda < \kappa} y_{\leq \lambda}.$$
The $z_\kappa$ are pairwise disjoint: if $\kappa_1 < \kappa_2$, then $$z_{\kappa_1} \leq y_{\leq \kappa_1} \perp z_{\kappa_2}.$$
Notice also that $\sum z_\kappa = 1$.  For if there were $z \in \p(\z(\M))$ with $z \perp (\sum z_\kappa)$, then let $\lambda$ be the least cardinal with $z y_{\leq \lambda} \neq 0$; by definition $z z_\lambda \neq 0$ as well, which contradicts the assumption.

We claim that
\begin{equation} \label{E:sup}
\sup_\alpha f^\alpha = \sum g_\kappa z_\kappa \triangleq f,
\end{equation}
where $g_\kappa = \kappa$ when $\kappa$ is infinite, and otherwise $g_\kappa = \sup_\alpha (g_\kappa^\alpha z_\kappa)$, which exists as the supremum of a bounded set in $\z(\M)_+$.

Next we show that $f \geq f^\alpha$ for any $\alpha$.  Fixing a cardinal $\lambda \leq \kappa_\M$,
\begin{equation} \label{E:cut}
z_\lambda f^\alpha = z_\lambda \left( \sum g_\kappa^\alpha z_\kappa^\alpha \right) = (z_\lambda y_{\leq \lambda}) \left( \sum g_\kappa^\alpha z_\kappa^\alpha \right) \leq z_\lambda \left( \sum_{\kappa \leq \lambda} g_\kappa^\alpha z_\kappa^\alpha \right).
\end{equation}
When $\lambda$ is infinite, we continue \eqref{E:cut} as
$$\leq \lambda z_\lambda = z_\lambda f.$$
When $\lambda$ is finite, we continue \eqref{E:cut} as
$$\leq z_\lambda (\lambda - 1) \left( \sum_{\kappa < \lambda}  z_\kappa^\alpha \right) + z_\lambda g_\lambda^\alpha z_\lambda^\alpha \leq z_\lambda g_\lambda = z_\lambda f.$$
Since $z_\lambda f^\alpha \leq z_\lambda f$ for all $\lambda$, $f \geq f^\alpha$.

Finally we check that if $h = \sum h_\kappa x_\kappa$ satisfies $h \geq f^\alpha$, $\forall \alpha$, then necessarily $h \geq f$.  Fixing a cardinal $\lambda \leq \kappa_\M$,
\begin{align*}
h_\lambda x_\lambda = x_\lambda h \geq x_\lambda f^\alpha, \: \forall \alpha \quad &\Rightarrow \quad x_\lambda \leq \sum_{\kappa \leq \lambda} z^\alpha_\kappa, \: \forall \alpha \\ &\Rightarrow \quad x_\lambda \leq \bigwedge_\alpha \left(\sum_{\kappa \leq \lambda} z^\alpha_\kappa \right) = y_{\leq \lambda}.
\end{align*}
This last inequality implies
\begin{equation} \label{E:cut2}
x_\lambda f = x_\lambda y_{\leq \lambda} f \leq x_\lambda \left( \sum_{\kappa \leq \lambda} g_\kappa z_\kappa \right).
\end{equation}
When $\lambda$ is infinite, we continue \eqref{E:cut2} as
$$\leq \lambda x_\lambda = x_\lambda h.$$
When $\lambda$ is finite, we continue \eqref{E:cut2} as
$$\leq x_\lambda (\lambda - 1) \left( \sum_{\kappa < \lambda} z_\kappa \right) + x_\lambda g_\lambda z_\lambda$$
and the inequality $h \geq f^\alpha$, $\forall \alpha$, allows us to compute further
$$ = x_\lambda (\lambda - 1) \left( \sum_{\kappa < \lambda} z_\kappa \right) + x_\lambda \left( \sup_\alpha g_\lambda^\alpha z_\lambda \right) \leq x_\lambda h.$$
Since $x_\lambda f \leq x_\lambda h$ for all $\lambda$, $f \leq h$.

This completes the proof that $f = \sup f^\alpha$.

As for the infimum of the $f^\alpha$, we first point out that we cannot write anything like
$$\bigwedge f^\alpha = 1 - \left(\bigvee (1-f^\alpha)\right),$$
which is a useful duality in $\p(\M)$.  There is no complementation in the lattices $\ps$ and $(\M_+/\approx)$, at least when $\M$ is not finite.  Instead we define
$$y_{\leq \kappa} = \bigvee_\alpha \left(\sum_{\lambda \leq \kappa} z_\lambda^\alpha \right)$$
and complete the rest of the proof similarly to the proof for the supremum.  (The substitute for \eqref{E:cut} should begin with ``$z_\lambda^\alpha f = \dots$"; for \eqref{E:cut2} should begin with ``$z_\lambda h = \dots$".) 
\end{proof}

\section{Application to representation theory} \label{S:rep}

In this section we reinterpret Theorem \ref{T:complat} in terms of the (normal Hilbert space) representations of $\M$.  Unless noted otherwise, we use ``isomorphism" in the sense of normed $\M$-modules, i.e.
$$\{\pi_1, \h_1\} \simeq \{\pi_2, \h_2\} \iff$$
$$\exists \text{ unitary }U: \h_1 \to \h_2: \qquad U\pi_1(x)U^* = \pi_2(x), \qquad \forall x \in \M.$$

It follows from the basic theory (see \cite[Sections 2.1-2]{JS} or \cite[Section 2]{S2003}) that any representation is (isomorphically) contained in a direct sum of copies of the standard form $\{\text{id}, L^2(\M)\}$.  We view $\oplus_I L^2(\M)$ as a row vector and think of the $\M$-action as multiplication on the left.  The commutant is right multiplication by $\B(\ell^2_I) \overline{\otimes} \M$, and the closed submodules are of the form $(\oplus_I L^2(\M)) q$, where $q \in \p(\B(\ell^2_I) \overline{\otimes} \M)$.  Two submodules are isomorphic if and only if the corresponding projections are equivalent.

This means that the isomorphism class of a representation corresponds to an equivalence class of projections in some amplification of $\M$.  Adding representations corresponds to adding equivalence classes.  As we have mentioned, the partial order can be defined in terms of the sum, so provided we make some kind of size restriction, we get an isomorphism of ordered monoids.  For example, if $\M$ is $\sigma$-finite, we obtain an identification between $(\p(\B(\ell^2) \overline{\otimes} \M)/\sim)$ and isomorphism classes of countably generated Hilbert $\M$-modules.  This all works for $L^p$ modules (\cite{JuS}), too, and is closely related to the $\K_0$ functor (\cite{Han,W-O}).

(Most of the ideas of the preceding two paragraphs were discussed by Breuer (\cite{B1968,B1969}), without making reference to standard forms.  He focused on the monoid generated by equivalence classes of finite projections, because the associated Grothendieck group, called the \textit{index group} of $\M$, is the natural carrier for the Fredholm theory of $\M$.  Olsen (\cite{O}) later combined Breuer's work with Tomiyama's dimension function to give a very general version of index theory in von Neumann algebras.) 

\begin{corollary}
Let $\{\pi_\alpha, \h_\alpha\}$ be a set of representations of a fixed von Neumann algebra $\M$.  Then there is a maximal representation of $\M$ which is (isomorphically) contained in all of these, and there is a minimal representation which (isomorphically) contains all of these.  Both are unique up to $\M$-module isomorphism.
\end{corollary}

\begin{proof}
Choose a large enough set $I$ so that for all $\alpha$, $\{\pi_\alpha, \h_\alpha\}$ is a subrepresentation of $\{\text{id}, \oplus_I L^2(\M)\}$.  The corollary follows from the preceding discussion and the fact that $(\p(\B(\ell^2_I) \overline{\otimes} \M)/\sim)$ is a complete lattice.
\end{proof}

In the early years of the subject, von Neumann algebras were generally given on Hilbert spaces, and the notion of $\M$-module isomorphism was therefore not in use.  Instead, one classified represented algebras up to the slightly weaker notion of \textit{spatial isomorphism}, which allows for an arbitrary isomorphism between the algebras.  (An $\M$-module isomorphism between representations $\{\pi_1, \h_1\}$ and $\{\pi_2, \h_2\}$ is a spatial isomorphism between von Neumann algebras $\{\pi_1(\M), \h_1\}$ and $\{\pi_2(\M), \h_2\}$ which induces the algebra isomorphism $\pi_2 \circ \pi_1^{-1}$.)  The question ``When is an algebraic isomorphism of represented von Neumann algebras spatial?", which is a noncommutative version of the fundamental problem of unitary equivalence for normal operators, is treated in detail in \cite{K1957}.  Also see \cite{Dig} for a projection-based approach to the existence of spatial isomorphisms.

Having said that, equivalence classes of representations/represented algebras were first studied by Murray and von Neumann (\cite[Chapter III]{MvN4}), using the coupling constant for finite factors.  The generalizations to coupling functions and arbitrary algebras were the motivations for the Griffin and Pallu de la Barri\`{e}re results featured in Section \ref{S:trace}.  The space-free approach was notably developed by the Japanese school of the 1950's.

Modulo spatial isomorphism, the set of equivalence classes of representations of a fixed von Neumann algebra may not even be partially ordered.  We mention the relevant example.  Let $\M$ be a type $\text{II}_\infty$ factor with dimension function $D$ and fundamental group $\Gamma \notin \{\{1\},(0,\infty)\}$.  (The existence of such an $\M$ remained in doubt until a breakthrough of Connes in 1980 (\cite{C}).  The fundamental group of a $\text{II}_\infty$ factor can be defined as
$$\{\lambda \in (0,\infty) \mid \exists \alpha \in \text{Aut}(\M), \; D \circ \alpha = \lambda D\},$$
with the group operation being multiplication.)  Kadison (\cite{K1955}) showed that for nonzero finite projections $p,q$, $L^2(\M)p$ is spatially isomorphic to $L^2(\M)q$ if and only if $\frac{D(p)}{D(q)} \in \Gamma$.  Since $\Gamma \neq (0,\infty)$, we may find nonzero finite projections $p,p'$ with $\frac{D(p)}{D(p')} \notin \Gamma$.  And $\Gamma \neq \{1\}$, so we may find spatial isomorphisms $L^2(\M)q_1 \simeq L^2(\M)p' \simeq L^2(\M)q_2$ with $q_1 \lneqq p \lneqq q_2$.  Therefore the spatial equivalence class of $L^2(\M)p$ both dominates and is dominated by that of $L^2(\M)p'$, yet the two are not equal.

At least for factors, this kind of pairing - $\text{II}_\infty$ algebra, $\text{II}_1$ commutant - is the only case where the two notions of equivalence differ.  Not coincidentally, the only choice required for $T$, $D$, and $\widehat{T}$ which cannot be standardized is the normalization on the finite elements in a $\text{II}_\infty$ summand.  (On a $\text{II}_\infty$ summand, one possible definition for ``normalization" is the inverse image of the identity, which is nothing but the equivalence class of the projection $p$ discussed in Remark \ref{R:trace}.)


\end{document}